\documentclass[12pt,a4paper,reqno]{amsart}
\pdfoutput=1
\usepackage[hidelinks,pdfstartview=FitH]{hyperref} 
\usepackage{fullpage}
\usepackage{amssymb,enumitem,mathrsfs,stmaryrd}
\usepackage[T1]{fontenc}
\usepackage{mathrsfs}
\usepackage[english]{babel}
\usepackage[justification=centering,skip=0pt]{caption}
\setenumerate{label=\textup{(\roman*)},itemsep=2pt,topsep=2pt,leftmargin=2em}
\setitemize{itemsep=2pt,topsep=2pt,leftmargin=2em}
\usepackage{tikz,pgfplots}
\usetikzlibrary{arrows,arrows.meta,calc,cd}

\tikzset{
>=Stealth,
node distance=2cm,
main node/.style={circle,inner sep=1pt},
freccia/.style={->,shorten >=2pt,shorten <=3pt},
ciclo/.style={out=130, in=50, loop, distance=2cm},
doppia/.style={-{Implies},double,double equal sign distance,shorten >=2pt,shorten <=2pt},
baseline=(current bounding box.center),
font=\small}

\usepackage{mathtools}
\mathtoolsset{showonlyrefs}

\allowdisplaybreaks[4]
\newcommand{\doi}[1]{\url{https://doi.org/#1}}
\newcommand{\isbn}[1]{\url{https://isbnsearch.org/isbn/#1}}

\newcommand{\web}[1]{\url{#1}}

\renewcommand{\emptyset}{\varnothing}
\numberwithin{equation}{section}

\newtheorem{prop}{Proposition}[section]
\newtheorem{thm}[prop]{Theorem}
\newtheorem{lemma}[prop]{Lemma}

\theoremstyle{remark}

\newtheorem{rem}[prop]{Remark}

\usepackage{etoolbox}
\AtEndEnvironment{rem}{\hfill\ensuremath{\square}}

\newcommand{\C}{\mathbb{C}}
\newcommand{\Z}{\mathbb{Z}}
\newcommand{\N}{\mathbb{N}}

\newcommand{\id}{\mathrm{id}}

\makeatletter
\def\@tocline#1#2#3#4#5#6#7{\relax
  \ifnum #1>\c@tocdepth 
  \else
    \par \addpenalty\@secpenalty\addvspace{#2}%
    \begingroup \hyphenpenalty\@M
    \@ifempty{#4}{%
      \@tempdima\csname r@tocindent\number#1\endcsname\relax
    }{%
      \@tempdima#4\relax
    }%
    \parindent\z@ \leftskip#3\relax \advance\leftskip\@tempdima\relax
    \rightskip\@pnumwidth plus4em \parfillskip-\@pnumwidth
    #5\leavevmode\hskip-\@tempdima
      \ifcase #1
       \or\or \hskip 1em \or \hskip 2em \else \hskip 3em \fi%
      #6 \hskip 0.5em \nobreak\relax
    \dotfill\hbox to\@pnumwidth{\@tocpagenum{#7}}\par
    \nobreak
    \endgroup
  \fi}
\makeatother

\title{On two groupoid models of quantum\\[5pt] complex projective spaces}

\author[de Castro]{Gilles G.~de Castro}
\author[D'Andrea]{Francesco D'Andrea}

\address[de Castro]{Departmento de Matem\'atica, Universidade Federal de Santa Catarina, Florian\'opolis, SC, 88040-900, Brazil}
\email{gilles.castro@ufsc.br}

\address[D'Andrea]{Dipartimento di Matematica e Applicazioni ``R.~Caccioppoli'',
Universit\`a di Napoli Federico II, and INFN -- Sezione di Napoli, Napoli, Italy}
\email{francesco.dandrea@unina.it}

\thanks{This research is part of the EU Staff Exchange project 101086394 ``Operator Algebras That One Can See''. We thank Jesse Reimann, Sophie E.~Zegers and Robert Yuncken for interesting discussions on this work.}

\begin{document}



\begin{abstract}
The quantum complex projective space $\C P^n_q$
is described by both the graph C*-algebra $C^*(\widetilde{F}_n)$ of an amplified graph $\widetilde{F}_n$, and the core $C^*(\widetilde{L}_{2n+1})^{U(1)}$ of the graph C*-algebra of a row-finite graph $\widetilde{L}_{2n+1}$ with no sinks. For $n=1$, there is a simple explicit formula for the isomorphism $C^*(\widetilde{F}_n)\to C^*(\widetilde{L}_{2n+1})^{U(1)}$ in terms of Cuntz-Krieger generators. For $n\geq 2$, there is an abstract proof of existence of such an isomorphism, but no explicit formula is known. In this short note, we provide an explicit formula for $n=2$.
Then, in the second part of the paper, we will use the groupoid model of these C*-algebras to give a geometric interpretation to this explicit formula and get some hints about its possible generalization.
\end{abstract}

\maketitle

\begin{center}
\begin{minipage}{0.8\textwidth}
\parskip=0pt\small\tableofcontents
\end{minipage}
\end{center}

\bigskip

\section{Introduction}\label{sec:1}

Quantum complex projective spaces associated with Vaksman--Soibelman quantum 
spheres~\cite{VS90} were first studied by Meyer in~\cite{Mey95}.
In 2002, Hong and Szyma{\'n}ski discovered that the C*-algebra of the $2n+1$ dimensional quantum sphere $S^{2n+1}_q$ is isomorphic
to a graph C*-algebra.
They provided two graphs, the graph $L_{2n+1}$ in Figure~\ref{fig:graphSnq} \cite[Theorem~4.4]{HS02} and the simpler graph $\widetilde{L}_{2n+1}$
in Figure~\ref{fig:sphereB} \cite[Prop.~A.1]{HS02}.
In the same paper, they argue that the core $C^*(L_{2n+1})^{U(1)}$,
that we interpret as describing the quantum complex projective plane $\C P_q^n$, is itself a graph C*-algebra \cite[Sect.~4.3]{HS02}, where the graph now is the one in Figure~\ref{fig:graphCPnq}, denoted by $F_n$.
For $n=1$, one can provide an explicit formula for the isomorphism $C^*(F_1)\to C^*(L_3)^{U(1)}$ (see e.g.~\cite[Prop.~4.3]{DZ26} in the special case $r=1$; notice also that $F_1=\widetilde{F}_1$ and $L_3=\widetilde{L}_3$). However, for $n\geq 2$ a similar formula proved to be quite elusive.
In \cite{DZ26} the result was generalized as follows:
given a finite directed acyclic graph $R$ without parallel edges, one can construct two directed graphs $E_R$ and $F_R$, one by adding a loop at every vertex of $R$ and the other by replacing every arrow of $R$ by countably infinitely many arrows. Then, $C^*(F_R)$ is isomorphic to the core of $C^*(E_R)$ \cite[Theorem~3.2]{DZ26}. This theorem is proved by using the property of these two C*-algebras of being AF, and by computing their dimension groups. In general, for the isomorphism $C^*(F_R)\to C^*(E_R)^{U(1)}$, no explicit formula in terms of Cuntz-Krieger generators is known.
As a special case, when $R$ is the graph with vertex set $\{1,\ldots,n+1\}$ and one edge from $i$ to $i+1$ for all $1\leq i\leq n$, then $E_R$ is the graph $\widetilde{L}_{2n+1}$, $F_R$ is the graph $\widetilde{F}_n$ in Figure~\ref{fig:CPnqB}, and we get the isomorphism
$C^*(\widetilde{F}_n)\cong C^*(\widetilde{L}_{2n+1})^{U(1)}$.

\begin{figure}[t]
\begin{tikzpicture}

\clip (-1.1,1.5) rectangle (11.2,-3);

\node[main node] (1) {};
\node[main node,right of=1] (2) {};
\node[main node,right of=2] (3) {};
\node[main node,right of=3] (4) {};
\node[right of=4] (5) {};
\node[main node,right of=5] (6) {};

\filldraw (1) circle (0.06) node[below left] {$1$};
\filldraw (2) circle (0.06);
\filldraw (3) circle (0.06);
\filldraw (4) circle (0.06);
\filldraw (6) circle (0.06) node[below right,inner sep=2pt] {$n+1$};

\path[freccia] (1) edge[ciclo] (1);
\path[freccia] (2) edge[ciclo] (2);
\path[freccia] (3) edge[ciclo] (3);
\path[freccia] (4) edge[ciclo] (4);
\path[freccia] (6) edge[ciclo] (6);

\path[freccia]
	(1) edge (2)
	(2) edge (3)
	(3) edge (4);
\path[freccia,dashed] (4) edge (6);
\path[freccia]
	(1) edge[bend right] (3)
	(1) edge[bend right=40] (4)
	(2) edge[bend right] (4);
\path[white]
	(1) edge[bend right=60] coordinate (7) (6)
	(2) edge[bend right=50] coordinate (8) (6)
	(3) edge[bend right=40] coordinate (9) (6);
\path[shorten <=1pt]
	(1) edge[out=-60,in=180] (7)
	(2) edge[out=-50,in=180] (8)
	(3) edge[out=-40,in=180] (9);
\path[->,dashed,shorten >=1pt]
	(7) edge[out=0,in=240] (6)
	(8) edge[out=0,in=230] (6)
	(9) edge[out=0,in=220] (6);
\path[freccia,dashed] (4) edge[bend right,dashed] (6);
\end{tikzpicture}

\caption{The graph $L_{2n+1}$.}
\label{fig:graphSnq}

\vspace{6mm}

\begin{tikzpicture}

      \node[main node] (1) {};
      \node[main node,right of=1] (2) {};
      \node[main node,right of=2] (3) {};
      \node[main node,right of=3] (4) {};
      \node[right of=4] (5) {};
      \node[main node,right of=5] (6) {};

      \filldraw (1) circle (0.06) node[below left] {$1$};
      \filldraw (2) circle (0.06);
      \filldraw (3) circle (0.06);
      \filldraw (4) circle (0.06);
      \filldraw (6) circle (0.06) node[below right] {$n+1$};

      \path[freccia] (1) edge[ciclo] (1);
      \path[freccia] (2) edge[ciclo] (2);
      \path[freccia] (3) edge[ciclo] (3);
      \path[freccia] (4) edge[ciclo] (4);
      \path[freccia] (6) edge[ciclo] (6);

      \path[freccia] (1) edge (2) (2) edge (3) (3) edge (4);
      \path[freccia,dashed] (4) edge (6);
\end{tikzpicture}

\caption{The graph $\widetilde{L}_{2n+1}$.}
\label{fig:sphereB}

\vspace{6mm}

\begin{tikzpicture}

\clip (-1.1,0.5) rectangle (11.2,-3);

\node (1) {};
\node (2) [right of=1] {};
\node (3) [right of=2] {};
\node (4) [right of=3] {};
\node (5) [right of=4] {};
\node (6) [right of=5] {};

\filldraw (1) circle (0.06) node[below left] {$1$};
\filldraw (2) circle (0.06);
\filldraw (3) circle (0.06);
\filldraw (4) circle (0.06);
\filldraw (6) circle (0.06) node[below right] {$n+1$};

\path[freccia]
	(1) edge node[fill=white] {$\infty$} (2)
	(2) edge node[fill=white] {$\infty$} (3)
	(3) edge node[fill=white] {$\infty$} (4);
\path[freccia,dashed] (4) edge (6);
\path[freccia]
	(1) edge[bend right] node[fill=white] {$\infty$} (3)
	(1) edge[bend right=40] node[fill=white] {$\infty$} (4)
	(2) edge[bend right] node[fill=white] {$\infty$} (4);
\path[white]
	(1) edge[bend right=60] coordinate (7) (6)
	(2) edge[bend right=50] coordinate (8) (6)
	(3) edge[bend right=40] coordinate (9) (6);
\path[shorten <=1pt]
	(1) edge[out=-60,in=180] node[fill=white] {$\infty$} (7)
	(2) edge[out=-50,in=180] node[fill=white] {$\infty$} (8)
	(3) edge[out=-40,in=180] node[fill=white] {$\infty$} (9);
\path[->,dashed,shorten >=1pt]
	(7) edge[out=0,in=240] (6)
	(8) edge[out=0,in=230] (6)
	(9) edge[out=0,in=220] (6);
\path[freccia,dashed] (4) edge[bend right,dashed] (6);
\end{tikzpicture}

\caption{The graph $F_n$.}
\label{fig:graphCPnq}

\vspace{6mm}

\begin{tikzpicture}

\clip (-1.1,0.5) rectangle (11.2,-0.5);

      \node (1) {};
      \node (2) [right of=1] {};
      \node (3) [right of=2] {};
      \node (4) [right of=3] {};
      \node (5) [right of=4] {};
      \node (6) [right of=5] {};

      \filldraw (1) circle (0.06) node[below left] {$1$};
      \filldraw (2) circle (0.06);
      \filldraw (3) circle (0.06);
      \filldraw (4) circle (0.06);
      \filldraw (6) circle (0.06) node[below right] {$n+1$};

      \path[freccia] (1) edge node[fill=white] {$\infty$} (2) (2) edge node[fill=white] {$\infty$} (3) (3) edge node[fill=white] {$\infty$} (4);
      \path[freccia,dashed] (4) edge (6);
\end{tikzpicture}

\caption{The graph $\widetilde{F}_n$.}
\label{fig:CPnqB}
\end{figure}

The aim of this paper is to construct an explicit isomorphism $C^*(\widetilde{F}_2)\to C^*(\widetilde{L}_5)^{U(1)}$, which is the first step for finding a formula for arbitrary $R$.
In the second part of the paper, we will use the groupoid model of these C*-algebras to give a geometric interpretation to the explicit formula for $n=2$, and get some hints for its generalization.
The busy reader can jump to Prop.~\ref{prop:isoA} and \ref{prop:isoB} for the explicit isomorphism of C*-algebras, or to Theorem~\ref{thm:41} for the isomorphism of topological groupoids.

The existence of an isomorphism between two topological groupoids (second countable, locally compact, Hausdorff, étale, essentially principal) is equivalent to the existence of a \emph{diagonal-preserving} isomorphism between the associated reduced groupoid \mbox{C*-algebras} \cite[Theorem~3.3]{CRSM21}. 
In the purely algebraic setting, an analogous result holds for Steinberg algebras (see \cite[Theorem~5.6]{Ste19} for the precise statement). One of the interesting aspects of such isomorphisms is that they establish a connection between the purely algebraic and the operator algebraic settings.

The plan of the paper is the following.
In Section~\ref{sec:2} we recall some basic definitions and results about graph C*-algebras.
Section~\ref{sec:3} is about the explicit isomorphism $C^*(\widetilde{F}_2)\to C^*(\widetilde{L}_5)^{U(1)}$.
In Section~\ref{sec:4} we recall the groupoid model of a graph C*-algebra.
In Section~\ref{sec:5} we present the explicit isomorphism between the groupoids underlying the C*-algebras $C^*(\widetilde{F}_2)$ and $C^*(\widetilde{L}_5)^{U(1)}$; this gives an independent proof of the isomorphism of C*-algebras above.
In Section~\ref{sec:6} we give an explicit formula for an isomorphism between groupoid in the general case: this is however only an isomorphism of the underlying algebraic groupoid; we will prove that it induces a homeomorphism between the spaces of objects, but it is not a homeomorphism on arrows.

\medskip

\noindent
{\bf Notations} Here $\N=\{0,1,2,\ldots\}$ is the set of natural numbers (including $0$). For $j,k\in\N$, $j\leq k$, we denote by $\binom{k}{j}:=\frac{k!}{j!\,(k-j)!}$ the binomial coefficient, with the convention that $0!:=1$. We set $\binom{k}{j}:=0$ if $j>k$. We adopt the convention that an empty sum is equal to $0$ and an empty product is equal to $1$.

\section{Graph C*-algebras}\label{sec:2}
We use the definition of graph C*-algebra in
\cite{FLR00}, in which sinks and infinite emitters are singular vertices. Let $E=(E^1,E^0,s,t)$ be a directed graph.
Recall that a vertex $v$ is a \emph{sink} if $s^{-1}(v)=\emptyset$,
it is an \emph{infinite emitter} if $|s^{-1}(v)|=\infty$, and it is a \emph{regular vertex} if $s^{-1}(v)$ is finite and non-empty. A vertex which is not regular is called \emph{singular}.
The graph C*-algebra $C^*(E)$ is the universal C*-algebras generated by a Cuntz-Krieger family $\{S_e,P_v\,|\,v\in E^0,e\in E^1\}$,
where $\{P_v\}$ is a family of orthogononal projections, $\{S_e\}$ a family of partial isometries with orthogonal ranges, and they satisfy the Cuntz-Krieger relations:
\begin{alignat}{2}
S_e^*S_e &=P_{t(e)} \qquad && \text{for all }e\in E^1, \label{eq:CK1} \\
S_eS_e^* &\leq P_{s(e)} && \text{for all }e\in E^1 , \label{eq:CK3} \\
\sum_{e\in E^1:\,s(e)=v}\!\! S_eS_e^* &=P_v && \text{for all regular }v\in E^0 . \label{eq:CK2}
\end{alignat}
The gauge action $\gamma:U(1)\to\mathrm{Aut}(C^*(E))$ is defined by $\gamma_z(S_e):=zS_e$ and $\gamma_z(P_v):=P_v$. The \emph{core} of $C^*(E)$ is the C*-subalgebra of all $U(1)$-invariant elements, denoted $C^*(E)^{U(1)}$.\noeqref{eq:CK2}

Notice that the conditions \eqref{eq:CK1} and that different partial isometries have orthogonal ranges can be combined into the single condition:\noeqref{eq:CK2p}
\begin{alignat}{2}
& S_e^*S_f =\delta_{e,f}P_{t(e)} \qquad && \text{for all }e,f\in E^1, \label{eq:CK1p} \\
\intertext{The conditions that $S_e$ is a partial isometry can then rephrased as:}
& S_eP_{t(e)} =S_e && \text{for all }e\in E^1 , \label{eq:CK2p} \\
\intertext{and \eqref{eq:CK3} can be rewritten as:}
& P_{s(e)}S_e =S_e && \text{for all }e\in E^1 . \label{eq:CK3p}
\end{alignat}

A \emph{finite path} in $E$ of length $n\geq 1$ is a sequence $\alpha=\alpha_1\alpha_2\ldots \alpha_n$ such that $\alpha_i\in E^1$ and $t(\alpha_i)=s(\alpha_{i+1})$ for all $i$.
A vertex is considered a path of length $0$.
Source and target of a path $\alpha=\alpha_1\alpha_2\ldots \alpha_n$ of length $n\geq 1$ are given by $s(\alpha):=s(\alpha_1)$ and $t(\alpha):=t(\alpha_n)$, respectively. If $\alpha$ is a vertex, we define $s(\alpha):=t(\alpha):=\alpha$.

We call $E^n$ the set of all paths in $E$ of length $n$, and $E^\star$ the set of all finite paths.

If $\alpha,\beta\in E^\star$ and $t(\alpha)=s(\beta)$, we denote by $\alpha\beta=\alpha_1\ldots\alpha_n\beta_1\beta_2\ldots$
their concatenation. If $\alpha=s(\beta)$ is a vertex, by definition $\alpha\beta:=\beta$. Similarly if $\beta=t(\alpha)$ is a vertex, then $\alpha\beta:=\alpha$.
By a \emph{cycle} we mean any finite path $\alpha$ such that $s(\alpha)=t(\alpha)$. 
A \emph{loop} is a cycle of length $1$.
If $\alpha$ is a cycle and $k\geq 1$, we define the cycle $\alpha^k$ in the obvious way, and adopt the convention that $\alpha^0=s(\alpha)$.

If $\alpha=\alpha_1\alpha_2\ldots\alpha_n$ is any path in $E$ of length $n\geq 1$,
we define as costumary
\[
S_\alpha:=S_{\alpha_1}S_{\alpha_2}\ldots S_{\alpha_n} ,
\]
and set $S_v:=P_v$ for $v\in E^0$ to uniform the notations.

\begin{rem}
If $E$ has no singular vertices, it follows from \eqref{eq:CK3} and induction on $n$ that
\begin{equation}\label{eq:mumu}
\sum_{\mu\in s^{-1}(v)\cap E^n} S_\mu S_\mu^*=P_v ,
\end{equation}
for all $v\in E$ and for all $n\geq 0$.
\end{rem}

\begin{rem}\label{rem:gen}
The core $C^*(E)^{U(1)}$ is spanned by products
$S_\alpha S_\beta^*$ with $t(\alpha)=t(\beta)=:v$ and $|\alpha|=|\beta|=:n\geq 0$
(see e.g.~\cite{R20}). If $v$ supports a loop $\ell$, then we can write
$S_\alpha S_\beta^*=(S_\alpha S_{\ell^n}^*) (S_\beta S_{\ell^n}^*)^*$.
If every vertex of $E$ supports a loop, then $C^*(E)^{U(1)}$ is generated
by elements of the form $S_\alpha S_{\ell^n}^*$, with $\alpha\in E^\star$,
$\ell$ a loop at $t(\alpha)$ and $n=|\alpha|$.
\end{rem}

\section{The C*-algebra of the quantum projective plane}\label{sec:3}

We label the edges and vertices of $\widetilde{L}_5$ and $\widetilde{F}_2$ as follows:
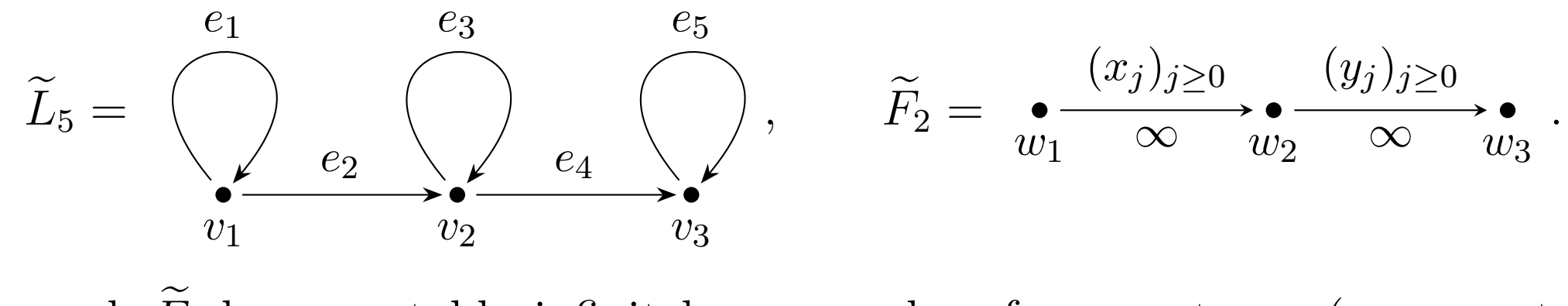
\begin{equation}\label{eq:shown}
\widetilde{L}_5=
\begin{tikzpicture}[font=\small,baseline=(current bounding box.center)]
\clip (-0.7,-0.5) rectangle (4.5,1.7);

\node[main node] (1) {};
\node (2) [main node,right of=1] {};
\node (3) [main node,right of=2] {};

\filldraw (1) circle (0.06) node[below=2pt] {$v_1$};
\filldraw (2) circle (0.06) node[below=2pt] {$v_2$};
\filldraw (3) circle (0.06) node[below=2pt] {$v_3$};

\path[freccia] (1) edge[ciclo] node[above] {$e_1$} (1);
\path[freccia] (2) edge[ciclo] node[above] {$e_3$} (2);
\path[freccia] (3) edge[ciclo] node[above] {$e_5$} (3);
\path[freccia] (1) edge node[above] {$e_2$} (2) (2) edge node[above] {$e_4$} (3);

\end{tikzpicture}
\;, \qquad
\widetilde{F}_2=
\begin{tikzpicture}[>=stealth,node distance=2cm,
main node/.style={circle,inner sep=2pt},
freccia/.style={->,shorten >=2pt, shorten <=2pt},
ciclo/.style={out=50, in=130, loop, distance=2cm},
font=\small,baseline={([yshift=-5pt]current bounding box.center)}]

      \node[main node] (1) {};
      \node (2) [main node,right of=1] {};
      \node (3) [main node,right of=2] {};

      \filldraw (1) circle (0.06) node[below=2pt] {$w_1$};
      \filldraw (2) circle (0.06) node[below=2pt] {$w_2$};
      \filldraw (3) circle (0.06) node[below=2pt] {$w_3$};

      \path[freccia] (1) edge node[below] {$\infty$} node[above] {$(x_j)_{j\geq 0}$} (2) (2) edge node[below] {$\infty$} node[above] {$(y_j)_{j\geq 0}$} (3);

\end{tikzpicture} .
\end{equation}
The graph $\widetilde{F}_2$ has countably infinitely many edges from $w_1$ to $w_2$ (resp.~$w_2$ to $w_3$) that are denoted by $x_j$ (resp.~$y_j$) with $j\in\N$.
To simplify the formulas, we will use the path algebra convention of denoting the Cuntz-Krieger generators by the same symbols used for the vertices and edges. That is, we write $v_i$ instead of $P_{v_i}$, $e_i$ instead of $S_{e_i}$, etc.

\begin{prop}\label{prop:isoA}
Fix a sequence $(N_j)_{j\geq 0}$ of natural numbers satisfying
\begin{equation}\label{eq:ass1}
N_j+j<N_{j+1}\;\forall\;j\in\N .
\end{equation}
Then, there exists a (unique) *-homomorphism $\Phi:C^*(\widetilde{F}_2)\to C^*(\widetilde{L}_5)^{U(1)}$ given on generators by:
\begin{gather}
  \Phi(w_i) =v_i , \qquad
  \Phi(y_j) = e_3^je_4(e_5^*)^{j+1},  \label{eq:elementsA} \\
  \Phi(x_j)  = e_1^je_2e_3{}^{N_j}(e_3^*)^{N_j+j+1}
    + \sum_{k=0}^{N_j+j} e_1^{j+1}e_2\,e_3^ke_4 (e_5^*)^{j+2}e_4^*(e_3^*)^k , \label{eq:elementsB}
\end{gather}
for all $i\in\{1,2,3\}$ and $j\in\N$.
This *-homomorphism is injective.
\end{prop}

\begin{proof}
To explain how the condition \eqref{eq:ass1} and the upper bound $N_j+j$ of the sum pop up, for the time being let us replace \eqref{eq:elementsB} with
\begin{equation}\label{eq:elementsC}
  \Phi(x_j)  = e_1^je_2e_3{}^{N_j}(e_3^*)^{N_j+j+1}
    + \sum_{k=0}^{M_j} e_1^{j+1}e_2\,e_3^ke_4 (e_5^*)^{j+2}e_4^*(e_3^*)^k , \tag{\ref{eq:elementsB}b}
\end{equation}
where $(M_j)_{j\geq 0}$ and $(N_j)_{j\geq 0}$ are two arbitrary sequences of natural numbers.

We start by noticing that the elements in \eqref{eq:elementsA} and
\eqref{eq:elementsC} belong to the core of $C^*(\widetilde{L}_5)$ and that the projections $\Phi(w_i)$, $i=1,2,3$, are orthogonal. Since every vertex in $\widetilde{F}_2$ is singular, in order to prove that the elements in \eqref{eq:elementsA} and \eqref{eq:elementsC} form a Cuntz-Krieger $\widetilde{F}_2$-family, it remains to check that the conditions \eqref{eq:CK1p}-\eqref{eq:CK3p} are satisfied.

We start with \eqref{eq:CK3p}.
Clearly $\Phi(w_1)\Phi(x_j)=\Phi(x_j)$
because $\Phi(x_j)$ is a sum of paths (in the extended graph) with source $v_1$ and $\Phi(w_1)=v_1$.
Similarly $\Phi(w_2)\Phi(y_j)=\Phi(y_j)$,  because $\Phi(y_j)$ has source $v_2$ and $\Phi(w_2)=v_2$.

Next, we pass to \eqref{eq:CK2p}. From $e_3^*v_2=e_3^*$,  $e_4^*v_2=e_4^*$ and $e_5^*v_3=e_5^*$ it follows that 
$\Phi(x_j)\Phi(w_2)=\Phi(x_j)$ and $\Phi(y_j)\Phi(w_3)=\Phi(y_j)$.

Finally, we check \eqref{eq:CK1p}. For all $j<k$,
\[
\Phi(y_j)^*\Phi(y_k)= e_5^{j+1}e_4^*e_3^{k-j}e_4(e_5^*)^{k+1}=0
\]
where the last equality follows from $e_4^*e_3=0$ and from $k-j\geq 1$. By adjunction, we deduce that
$\Phi(y_j)^*\Phi(y_k)=0$ for $j>k$ as well. On the other hand, if $j=k$ we get
\begin{align*}
\Phi(y_j)^*\Phi(y_j) &= e_5^{j+1}e_4^*v_2e_4(e_5^*)^{j+1}=e_5^{j+1}e_4^*e_4(e_5^*)^{j+1} \\
&=e_5^{j+1}v_3(e_5^*)^{j+1}=
e_5^{j+1}(e_5^*)^{j+1}=w_3=\Phi(v_3) .
\end{align*}
Next, $\Phi(x_j)$ and $\Phi(y_k)$ have orthogonal range projections $v_1$ and $v_2$,
hence $\Phi(x_j)^*\Phi(y_k)=0$ and $\Phi(y_k)^*\Phi(x_j)=0$.
Finally, call
\begin{equation}\label{eq:AB}
A_j:=e_1^je_2\,e_3^{N_j}(e_3^*)^{N_j+j+1} , \qquad\quad B_{i,j}:=e_1^{j+1}e_2\,e_3^ie_4 (e_5^*)^{j+2}e_4^*(e_3^*)^i ,
\end{equation}
so that \[ \Phi(x_j)  =A_j+\sum_{i=0}^{M_j}B_{i,j} . \]
Then, for all $j<k$, one has
\begin{align*}
A_j^*A_k &=e_3^{N_j+j+1}(e_3^*)^{N_j} [e_2^* e_1^{k-j}] e_2 e_3^{N_k}(e_3^*)^{N_k+k+1} =0 ,
\\
A_j^*B_{i',k} &=e_3^{N_j+j+1}(e_3^*)^{N_j} [e_2^* e_1^{k-j+1}] e_2e_3^{i'}e_4 (e_5^*)^{k+2}e_4^*(e_3^*)^{i'} =0 ,
\\
B_{i,j}^*B_{i',k} &=e_3^ie_4e_5^{j+2}e_4^* (e_3^*)^i [e_2^* e_1^{k-j}] e_2e_3^{i'}e_4 (e_5^*)^{k+2}e_4^*(e_3^*)^{i'} =0 .
\end{align*}
These are all $0$ because the square brackets contain a positive power of $e_1$, and $e_2^*e_1=0$. 
For the same reason, if $j<k-1$, then
\[
B_{i,j}^*A_k=e_3^ie_4e_5^{j+2}e_4^* (e_3^*)^i[e_2^*e_1^{k-j-1}]e_2e_3^{N_k}(e_3^*)^{N_k+k+1} =0 .
\]
On the other hand, if $j=k-1$,
\begin{align*}
B_{i,j}^*A_{j+1} &=e_3^ie_4e_5^{j+2}e_4^* (e_3^*)^ie_2^*e_2e_3^{N_{j+1}}(e_3^*)^{N_{j+1}+j+2}   \\
&=e_3^ie_4e_5^{j+2}e_4^* (e_3^*)^ie_3^{N_{j+1}}(e_3^*)^{N_{j+1}+j+2} .
\end{align*}
If $i\geq N_{j+1}$ this does not vanish. If $i<N_{j+1}$, on the other hand, 
the last line is $0$ because $e_4^*e_3=0$.
We need to assume that $i<N_{j+1}$ for all $0\leq i\leq M_j$,
that is
\begin{equation}\label{eq:ass1}
M_j<N_{j+1}\;\forall\;j\in\N .
\end{equation}
Under this assumption, we get $\Phi(x_j)^*\Phi(x_k)=0$ for all $j<k$, and by adjunction for all $j>k$ as well.
Finally, when $j=k$,
\[
B_{i,j}^*A_j =e_3^ie_4e_5^{j+2}e_4^* (e_3^*)^ie_2^*
e_1^*e_2\,e_3^{N_j}(e_3^*)^{N_j+j+1}=0
\]
because $e_1^*e_2=0$. By adjunction, $A_j^*B_{i,j}=0$. Moreover,
\[
B_{i,j}^*B_{i',j}=e_3^ie_4e_5^{j+2}e_4^* (e_3^*)^ie_3^{i'}e_4(e_5^*)^{j+2}e_4^*(e_3^*)^{i'} = 0
\]
if $i\neq i'$, either because $e_4^*e_3=0$ ($i'>i$) or $e_3^*e_4=0$ ($i>i'$).
Hence,
\[
\Phi(x_j)^*\Phi(x_j)=A_j^*A_j+\sum_{i=0}^{M_j} B_{i,j}^*B_{i,j}.
\]
An explicit computation gives
\begin{align*}
A_j^*A_j&=e_3^{N_j+j+1}(e_3^*)^{N_j+j+1} ,
\\
B_{i,j}^*B_{i,j}&=e_3^ie_4e_5^{j+2}(e_5^*)^{j+2}e_4^*(e_3^*)^i .
\end{align*}
Note that
\[
B_{i,j}^*B_{i,j}=e_3^ie_4e_5^{N_j+j-i}(e_5^*)^{N_j+j-i}e_4^*(e_3^*)^i
\]
for all $i\leq N_j+j$ (since $e_5e_5^*=e_5^*e_5=v_3$).
If $M_j\leq N_j+j$, then
\begin{equation}\label{eq:summu}
\Phi(x_j)^*\Phi(x_j)=\sum_{\mu} \mu\mu^* ,
\end{equation}
where the sum is over (some) paths $\mu$ of length $N_j+j+1$ with source $v_2$. If
\begin{equation}\label{eq:ass2}
M_j=N_j+j ,
\end{equation}
in the sum \eqref{eq:summu} all such paths appear:
indeed, either $\mu=e_3^{N_j+j+1}$
if the target is also $v_2$, or $\mu$ does $i$ loops around $e_3$, with
$0\leq i\leq M_j=N_j+j$, then moves to $v_3$ and does $N_j+j-i$ loops around $e_5$. From \eqref{eq:mumu}
it follows that $\Phi(x_j)^*\Phi(x_j)=v_2=\Phi(w_2)$.

This proves that $\Phi$ is well-defined. Since $\widetilde{F}_2$ has no cycles, the condition (L) is satisfied (every cycle has an exit), and injectivity of $\Phi$ follows from the Cuntz-Krieger Uniqueness Theorem \cite[Theorem 2]{FLR00}.
\end{proof}

\begin{lemma}\label{lemma:Phixy}
The *-homomorphism in Prop.~\ref{prop:isoA} is given, on paths of length $2$, by the formula:
\begin{equation}\label{eq:Phijk}
\Phi(x_jy_k)  =
\begin{cases}
 e_1^{j+1}e_2\,e_3^ke_4 (e_5^*)^{j+k+3} & \text{if } k<N_j+j+1 , \\[2pt]
 e_1^je_2\,e_3^{k-j-1}e_4(e_5^*)^{k+1} & \text{if } k\geq N_j+j+1 ,
 \end{cases}
\end{equation}
for all $j,k\in\N$.
\end{lemma}

\begin{proof}
Let us decompose again $\Phi(x_j)  =A_j+\sum_{i=0}^{N_j+j}B_{i,j}$,
with $A_j$ and $B_{i,j}$ as in \eqref{eq:AB}.
If  $N_j+j+1>k$, then
$A_j\Phi(y_k)=0$ because $e_3^*e_4=0$. If $N_j+j+1\leq k$, then
\[
A_j\Phi(y_k)=e_1^je_2\,e_3^{k-j-1} e_4(e_5^*)^{k+1} .
\]
On the other hand, from $e_4^*(e_3^*)^ie_3^ke_4=\delta_{i,k}v_3$
we deduce that
\[
\sum_{i=0}^{N_j+j}B_{i,j}\Phi(y_k)
= e_1^{j+1}e_2\,e_3^ke_4 (e_5^*)^{j+k+3}
\]
if $k<N_j+j+1$, and the product is $0$ otherwise.
\end{proof}

\begin{lemma}\label{lemma:sigma}
The map
\[
\sigma:\N^2\to\N^2 , \qquad
\sigma(j,k) :=\begin{cases}
 (j+1,k) & \text{if } k<N_j+j+1 , \\[2pt]
 (j,k-j-1) & \text{if } k\geq N_j+j+1 ,
\end{cases}
\]
is a bijection if and only if
\begin{equation}\label{eq:Nj}
N_j=\binom{j+1}{2} \;\forall\;j\in\N .
\end{equation}
In the latter case, its inverse is given by
\[
\sigma^{-1}:\N^2\to\N^2, \qquad
\sigma^{-1}(r,s) :=\begin{cases}
 (r-1,s) & \text{if } s<N_r , \\[2pt]
 (r,s+r+1) & \text{if } s\geq N_r .
\end{cases}
\]
\end{lemma}

\begin{proof}
The map $\sigma$ is well-defined, since $k-j-1\geq 0$ when $k\geq N_j+j+1$. 

We prove the implication ``$\Leftarrow$'' first. Assume that \eqref{eq:Nj} holds.
The map $\sigma^{-1}$ is then well-defined, since $\binom{r+1}{2}>s\geq 0$ implies that $r-1\geq 0$.

If $k<\binom{j+2}{2}$, then $r:=j+1$ and $s:=k$ satisfy
\[
s=k<\binom{j+2}{2}=\binom{r+1}{2} .
\]
In this case, when computing $\sigma^{-1}(\sigma(j,k))$ we use the first branch of both $\sigma$ and $\sigma^{-1}$ and easily verify that $\sigma^{-1}(\sigma(j,k))=(j,k)$.
If $k\geq\binom{j+2}{2}$, then $r:=j$ and $s:=k-j-1$ satisfy
\[
s=k-j-1\geq\binom{j+2}{2}-\binom{j+1}{1}=\binom{j+1}{2}=\binom{r+1}{2} .
\]
In this case, we compute $\sigma^{-1}(\sigma(j,k))$ using the second branch of both $\sigma$ and $\sigma^{-1}$ finding again $\sigma^{-1}(\sigma(j,k))=(j,k)$.
Thus, $\sigma^{-1}\circ\sigma=\id_{\N^2}$. Similarly one proves $\sigma\circ\sigma^{-1}=\id_{\N^2}$.

Finally, we prove the implication ``$\Rightarrow$''.
For $i\in\N$ let us see which points $(i,m)$ are in the image of $\sigma$. If $i=0$, then we get all points $(0,m)$ with $m\geq N_0$.
If $\sigma$ is surjective, it must be then
\begin{equation}\label{eq:N0}
N_0=0 .
\end{equation}
For $i\geq 1$, in the image we find points $(i,m)$ with either $m\geq N_i$ or $m<N_{i-1}+i$. In order to get all $m\in\N$, it must be
\begin{equation}\label{eq:Ni}
N_i=N_{i-1}+i+1 \;\forall\;i\geq 1.
\end{equation}
From \eqref{eq:N0} and the recursive equation \eqref{eq:Ni} we see that $N_j$ is the binomial \eqref{eq:Nj}.
\end{proof}

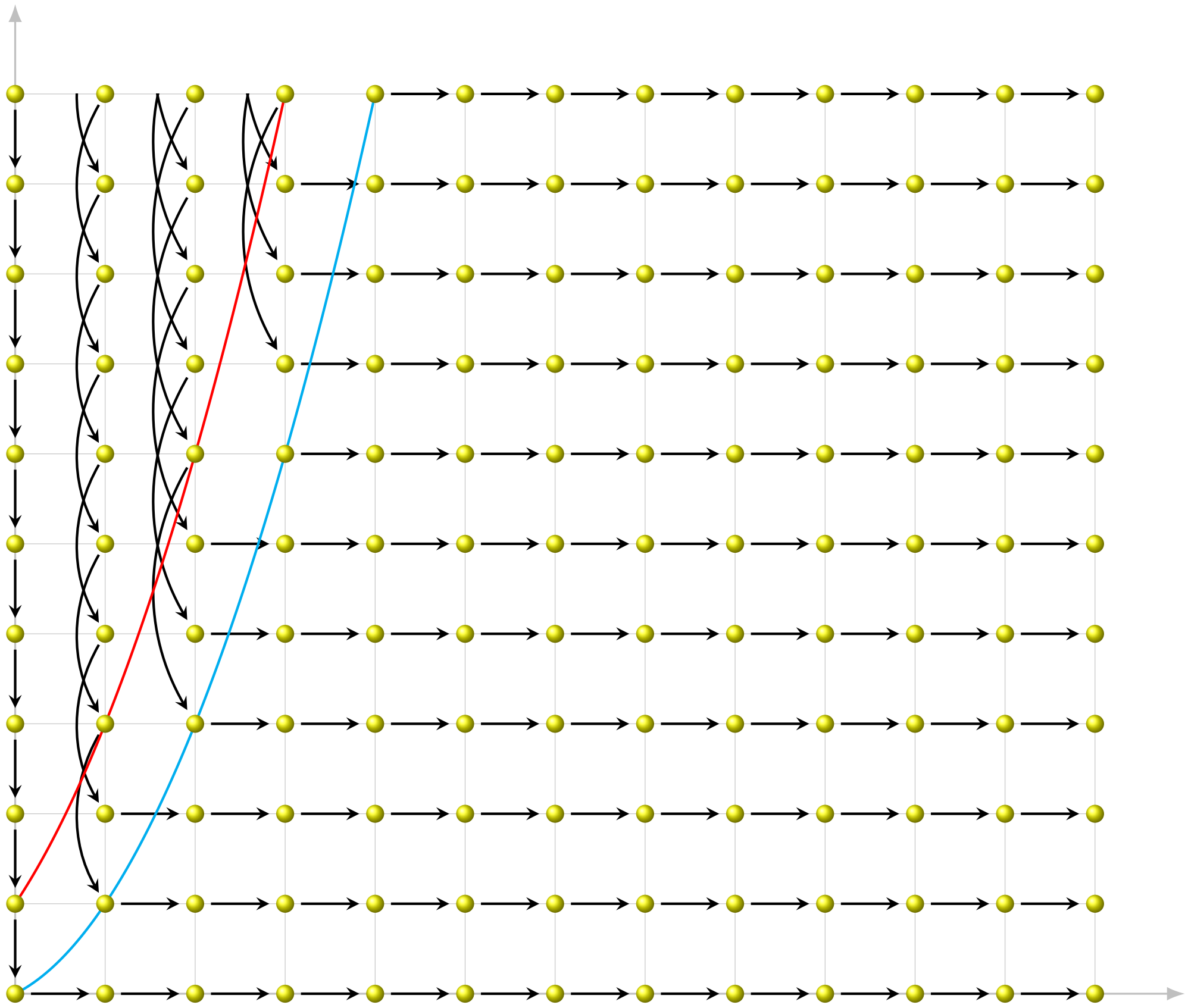
\begin{figure}[t]
\begin{tikzpicture}

\draw[very thin,color=gray!50] (0,0) grid (12,10);
\draw[gray!50,-Latex,semithick] (0,0) -- ++(13,0);
\draw[gray!50,-Latex,semithick] (0,0) -- ++(0,11);

\foreach \x in {4,5,...,11} \foreach \y in {0,1,...,10}
	\draw[thick,-stealth,shorten >=5pt, shorten <=5pt] (\x,\y) -- ++ (1,0);

\foreach \x in {3} \foreach \y in {0,1,...,9}
	\draw[thick,-stealth,shorten >=5pt, shorten <=5pt] (\x,\y) -- ++ (1,0);

\foreach \x in {2} \foreach \y in {0,1,...,5}
	\draw[thick,-stealth,shorten >=5pt, shorten <=5pt] (\x,\y) -- ++ (1,0);

\foreach \x in {1} \foreach \y in {0,1,...,2}
	\draw[thick,-stealth,shorten >=5pt, shorten <=5pt] (\x,\y) -- ++ (1,0);

\foreach \x in {0} \foreach \y in {0}
	\draw[thick,-stealth,shorten >=5pt, shorten <=5pt] (\x,\y) -- ++ (1,0);

\foreach \x in {0} \foreach \y in {1,2,...,10}
	\draw[thick,-stealth,shorten >=5pt, shorten <=5pt] (\x,\y) -- ++ (0,-1);

\begin{scope}
\clip (0,0) rectangle (3.5,10);

\foreach \x in {1} \foreach \y in {3,4,...,12}
	\path[thick,-stealth,shorten >=4pt, shorten <=4pt] (\x,\y) edge[bend right=30] ++ (0,-2);

\foreach \x in {2} \foreach \y in {6,7,...,13}
	\path[thick,-stealth,shorten >=5pt, shorten <=5pt] (\x,\y) edge[bend right=30] ++ (0,-3);

\foreach \x in {3} \foreach \y in {10,11,...,14}
	\path[thick,-stealth,shorten >=5pt, shorten <=5pt] (\x,\y) edge[bend right=30] ++ (0,-3);
\end{scope}

\draw[thick,color=red, domain=0:3, samples=200] plot (\x, {(\x + 1)*(\x + 2)/2});
\draw[thick,color=cyan, domain=0:4, samples=200] plot (\x, {\x*(\x + 1)/2});

\foreach \x in {0,1,...,12} \foreach \y in {0,1,...,10}
	\shade[ball color=yellow] (\x,\y) circle (0.1);

\end{tikzpicture}

\bigskip

\caption{The permutation $\sigma$.}\label{fig:5}
\end{figure}

The map $\sigma$ is illustrated in Figure \ref{fig:5}. The red curve is the parabola $k=(j+1)(j+2)/2$, the red one is the parabola $k=j(j+1)/2$. All the pairs $(j,k)\in\N^2$ strictly below the red parabola are shifted to the right by $1$. All the pairs $(j,k)$ above or on the red parabola are shifted down of $j+1$ slots.

\begin{prop}\label{prop:isoB}
If $N_j$ is as in \eqref{eq:Nj}, then the
*-homomorphism $\Phi:C^*(\widetilde{F}_2)\to C^*(\widetilde{L}_5)^{U(1)}$ in Prop.~\ref{prop:isoA} is an isomorphism.
\end{prop}

\begin{proof}
Let $N_j$ be given by \eqref{eq:Nj}.
We must prove the surjectivity of $\Phi$.
Since $\widetilde{L}_5$ satisfies the hypothesis of Remark \ref{rem:gen}, the core $C^*(\widetilde{L}_5)^{U(1)}$
is generated by elements $\mu(\ell^*)^{|\mu|}$ with $\mu$ a finite path and $\ell$ a loop based at $t(\mu)$.
It is enough to show that these elements are in the image of $\Phi$.
When $\mu$ is a vertex, the statement reduces to $\Phi(w_i) =v_i$, $i=1,2,3$.
Here is the list of all the remaining elements that we need to find in the image ($p,q,r\geq 0$):
\begin{align}
& e_1^{p+1}(e_1^*)^{p+1} & \text{(source $v_1$ and target $v_1$)} , \label{eq:el1} \\
& e_1^pe_2e_3^q(e_3^*)^{p+q+1} & \text{(source $v_1$ and target $v_2$)} , \label{eq:el2} \\
& e_1^pe_2e_3^qe_4e_5^r(e_5^*)^{p+q+r+2} 
=e_1^pe_2e_3^qe_4(e_5^*)^{p+q+2} 
& \text{(source $v_1$ and target $v_3$)} , \label{eq:el3}  \\
& e_3^{p+1}(e_3^*)^{p+1} & \text{(source $v_2$ and target $v_2$)} , \label{eq:el4} \\
& e_3^pe_4e_5^q(e_5^*)^{p+q}=
e_3^pe_4(e_5^*)^p
 & \text{(source $v_2$ and target $v_3$)} , \label{eq:el5} \\
& e_5^{p+1}(e_5^*)^{p+1}=v_3 & \text{(source $v_3$ and target $v_3$)} . \label{eq:el6} 
\end{align}
The element \eqref{eq:el6} is the image of $w_3$.
The element \eqref{eq:el5} is the image of $y_p$.
From \eqref{eq:mumu}, for all $n\geq 0$:
\begin{equation}\label{eq:33}
e_3^n(e_3^*)^n = v_2-\sum\nolimits_{\nu}\nu\nu^*=
v_2-\sum\nolimits_{\nu}\big(\nu(e_3^*)^n\big)\big(\nu(e_3^*)^n\big)^* ,
\end{equation}
where the sum is over all paths $\nu$ with source $v_2$ of length $n$ and different from $e_3^n$. Since $\nu(e_3^*)^n$ is in the image of $\Phi$, see \eqref{eq:el5}, then
\eqref{eq:el4} is in the image as well. A similar argument using \eqref{eq:mumu}
shows that, if \eqref{eq:el2} and \eqref{eq:el3} are in the image of $\Phi$, then so is
\eqref{eq:el1}.

From Lemma \ref{lemma:Phixy} and Lemma \ref{lemma:sigma} we see that
\[
\{\Phi(x_jy_k)\,|\,j,k\in\N\}=\{  e_1^pe_2e_3^qe_4 (e_5^*)^{p+q+2} | p,q\in\N \} ,
\]
which takes care of \eqref{eq:el3}. We are left with \eqref{eq:el2}.

Let $\mu_{p,q}=e_1^pe_2e_3^q(e_3^*)^{p+q+1}$ be the element in \eqref{eq:el2}.
From \eqref{eq:CK2}, $e_3e_3^*=v_2-e_4e_4^*$. Hence, 
for all $q\geq 1$:
\[
\mu_{p,q}=\mu_{p,q-1}-
e_1^pe_2e_3^{q-1}e_4e_4^*(e_3^*)^{p+q}  .
\]
The second summand can be rewritten as $(e_1^pe_2e_3^{q-1}e_4(e_5^*)^{p+q+1})(e_3^{p+q}e_4(e_5^*)^{p+q+1})^*$, that is a product of an element of the form \eqref{eq:el3} with the adjoint of one of the form \eqref{eq:el5}. From this, we deduce that, for all $p\geq 0$ and $q\geq 1$,
\begin{equation}\label{eq:iff}
\mu_{p,q}\in\mathrm{im}(\Phi) \iff \mu_{p,q-1}\in\mathrm{im}(\Phi) .
\end{equation}
Next, let us write $\Phi(x_j)  =A_j+\sum B_{i,j}$,
with $A_j$ and $B_{i,j}$ as in \eqref{eq:AB}. Since
\[
B_{i,j}=\big(e_1^{j+1}e_2\,e_3^ie_4 (e_5^*)^{i+j+3}\big) \big(e_3^ie_4(e_5^*)^{i+1}\big)^*
\]
is a product of an element of the form \eqref{eq:el3} with the adjoint of an element of the form \eqref{eq:el5}, it belongs to the image of $\Phi$. But then
\[
\sum B_{i,j}-\Phi(x_j)=A_j=\mu_{j,\binom{j+1}{2}}
\]
belongs to $\Phi$ as well. From the implication ``$\Rightarrow$'' in \eqref{eq:iff} it follows that $\mu_{p,q}\in\mathrm{im}(\Phi)$ for all $q\leq\binom{p+1}{2}$, and
the implication ``$\Leftarrow$'' in \eqref{eq:iff} it follows that $\mu_{p,q}\in\mathrm{im}(\Phi)$ for all $q>\binom{p+1}{2}$.
\end{proof}

\section{The groupoid model of a graph C*-algebra}\label{sec:4}
We assume that the reader is familiar with étale groupoids (see e.g.~\cite{SSW20}).
Here we will only fix the notations and recall the definition and some basic properties of the boundary path groupoid $\mathcal{G}_E$ of a directed graph $E$.

\subsection{The path space of a graph}\label{sec:2.2}
A \emph{right-infinite path} in a directed graph $E$ is right-infinite sequence $\alpha=\alpha_1\alpha_2\ldots \alpha_n\ldots$ with $\alpha_i\in E^1$ and $t(\alpha_i)=s(\alpha_{i+1})$ for all $i$. Its source is $s(\alpha):=s(\alpha_1)$, and
the concatenation $\beta\alpha$ makes sense whenever $\beta$ is a finite path
and $t(\beta)=s(\alpha)$.
If $\ell$ is a loop or a cycle, we define the right-infinite path $\ell^\infty$ in the obvious way.
We call $E^\infty$ the set of all right-infinite paths, $E^\star_{\mathrm{sing}}$ the set of finite paths that end at a singular vertex, and
\[
\partial E=E^\infty\cup E^\star_{\mathrm{sing}}
\]
the set of \emph{boundary paths}.
If $\alpha$ is a finite path, we define the \emph{cylinder set} $C_E(\alpha)$ as
\[
C_E(\alpha):=\big\{\alpha x\,\big|\,x\in \partial E,s(x)=t(\alpha)\big\} .
\]
If $S$ is a finite set of edges with source $t(\alpha)$, we define the corresponding \emph{generalised cylinder set} $C_E(\alpha,S)$ as
\[
C_E(\alpha,S):=C_E(\alpha)\setminus\bigcup_{e\in S}C_E(\alpha e) .
\]
The generalized cylinder sets form a basis of compact open sets of a totally disconnected locally compact Hausdorff topology on $\partial E$. With this topology, we call $\partial E$ the \emph{boundary path space} of $E$.
If $E^0$ is finite, then $\partial E$ is compact.
If $E^1$ is countable, then $\partial E$ is second-countable. See \cite[Sect.~2.3.2]{R20}.

Notice that, if $\alpha\in E^\star$ and $t(\alpha)$ is a regular vertex, then $\alpha$ is not a boundary path and $C_E(\alpha,S)=\bigcup_{e\notin S}C_E(\alpha e)$.
These sets are redundant, since they are unions of cylinder sets. We have then the following remark:

\begin{rem}\label{rem:similarly}
\begin{enumerate}
\item
A basis of the topology of $\partial E$ is given by the cylinder sets $C_E(\alpha)$, for all $\alpha\in E^\star$, and by the generalized cylinder sets
\[
C_E(\alpha,S)=\{\alpha\}\cup\bigcup_{e\in S^c}C_E(\alpha e)
\]
for all $\alpha$ whose target is an infinite emitter and for all cofinite subsets $S^c$ of $s^{-1}(t(\alpha))$.

\item
If $t(\alpha)$ is a sink, then $\{\alpha\}=C_E(\alpha)$ is open.\vspace*{-14pt}
\end{enumerate}
\end{rem}

\subsection{The boundary path groupoid}
We now recall the construction of the boundary path groupoid $\mathcal{G}_E$ of a graph $E$.

Let $k\in\Z$. Two paths $x,y\in\partial E$ are \emph{tail equivalent with lag $k$}, written $x\sim_ky$, if there exists $\alpha,\beta\in E^\star$ and $z\in\partial E$ such that
\[
t(\alpha)=t(\beta)=s(z) , \quad
|\alpha|-|\beta|=k , \quad
x=\alpha z , \quad
y=\beta z .
\]
In general, the lag of a pair of infinite paths is not unique. On the other hand, for tail equivalent paths $x=\alpha z$ and $y=\beta z$ ending at a singular vertex, the lag $k=|\alpha|-|\beta|=|x|-|y|$ is uniquely determined by $x$ and $y$, and two  paths ending at a singular vertex are tail equivalent if and only if they have the same target.

The \emph{boundary path groupoid} of $E$ has arrows
\[
\mathcal{G}_E^1:=\big\{(x,k,y)\,\big|\,x,y\in\partial E,x\sim_ky\big\} .
\]
The composition of $(x,k,y)$ with $(x',k',y')$ is defined when the \emph{domain} $y$ of the first morphism is equal to the \emph{codomain} $x'$ of the second morphism, and it is given by
\[
(x,k,y)(y,k',y') :=(x,k+k',y') .
\]
We will identify the unit space $\mathcal{G}_E^0:=\big\{(x,0,x)\,\big|\,x\in\partial E\big\}$ with $\partial E$.

The topology is defined as follows. For every pair of finite paths $(\alpha,\beta)$ with the same target, we define
\[
Z_E(\alpha,\beta):=\big\{(\alpha z,|\alpha|-|\beta|,\beta z)\,\big|\,z\in\partial E,s(z)=t(\alpha)=t(\beta)\big\} ,
\]
and for $S$ a finite set of edges with source $t(\alpha)$ we define
\[
Z_E(\alpha,\beta,S):=Z_E(\alpha,\beta)\setminus\bigcup_{e\in S}Z_E(\alpha e,\beta e).
\]
The sets $Z_E(\alpha,\beta,S)$ are the basis of a topology 
on $\mathcal{G}_E^1$ inducing on $\mathcal{G}^0_E$
the same topology as in Section \ref{sec:2.2}.
With this topology, $\mathcal{G}_E$ is a Hausdorff ample groupoid and $\{Z_E(\alpha,\beta,S)\}$ is a basis of compact open bisections \cite[Theorem 2.4]{R20}.
Similarly to Remark \ref{rem:similarly}, one has:

\begin{rem}\label{rem:resim}
\begin{enumerate}
\item\label{rem:resimA}
A basis of the topology of $\mathcal{G}_E^1$ is given by the sets $Z_E(\alpha,\beta)$, for all $\alpha,\beta\in E^\star$ with the same target, and by the sets
\[
Z_E(\alpha,\beta,S)=\{(\alpha,\beta)\}\cup\bigcup_{e\in S^c}Z_E(\alpha e,\beta e)
\]
for all $\alpha,\beta\in E^\star$ such that $t(\alpha)=t(\beta)=:v$ is an infinite emitter and for all cofinite subsets $S^c$ of $s^{-1}(v)$.

\item\label{rem:resimB}
If $t(\alpha)=t(\beta)$ is a sink, then $\{(\alpha,\beta)\}=Z_E(\alpha,\beta)$ is open.

\item\label{rem:resimC}
If $\alpha,\beta\in F^\star$ and $t(\alpha)=t(\beta)=:v$ is an infinite emitter, one can easily verify that, for every finite set $S$ of edges with source $v$, one has
\[
Z_E(\alpha,\beta,S)=Z_E(\alpha,v,S)\cdot Z_E(\beta,v,S)^{-1} . \vspace*{-20pt}
\]

\end{enumerate}
\end{rem}

\subsection{The convolution algebra}

The groupoid model of a graph C*-algebra $C^*(E)$ is discussed in \cite{KPRR97} when $E$ is locally finite with no sinks and no sources, in \cite{P02} when $E$ is any graph without sinks, and in \cite{BCM20} for an arbitrary graph. A similar result for Leavitt path algebras can be found in \cite{R20}, where it is proved that the Steinberg algebra of $\mathcal{G}_E$ is isomorphic to the Leavitt path algebra $L_{\C}(E)$. By passing to the C*-enveloping algebras, this gives an independent proof that the full groupoid C*-algebra $C^*(\mathcal{G}_E)$ is isomorphic to the graph C*-algebra $C^*(E)$.
If $E$ is countable and without sinks, then $\mathcal{G}_E$ is amenable \cite[Theorem 4.2]{P02}, which means that the reduced and full groupoid \mbox{C*-algebras} are isomorphic.

\begin{thm}[see e.g.~\protect{\cite[Theorem~2.7]{R20}}]\label{thm:isoEG}
Let $E$ be a countable graph and $\mathcal{G}_E$ its boundary path groupoid. Then $C^*(\mathcal{G}_E)\cong C^*(E)$, with isomorphism $\Psi_E$ defined by
\[
\Psi_E(P_v)=1_{Z_E(v,v)} ,
\qquad
\Psi_E(S_e)=1_{Z_E(e,t(e))} ,
\]
where $\{P_v,S_e\}$ is the Cuntz--Krieger $E$-family of $C^*(E)$ 
and $1_X$ denotes the characteristic function of $X\subseteq\mathcal{G}_E^1$.
\end{thm}

By induction on the length of the paths one easily shows that
\begin{equation}\label{eq:SZ}
\Psi_E(S_\alpha S_\beta^*)=1_{Z_E(\alpha,\beta)} ,
\end{equation}
for all $\alpha,\beta\in E^\star$ with $t(\alpha)=t(\beta)$.

The map $\Psi_E$ transforms the gauge action into a strongly $U(1)$-continuous action $\gamma$ given on $f\in C_c(\mathcal{G}_E^1)$ by
\[
(\gamma_zf)(g)=z^{c(z)}f(g) ,\qquad\quad\forall z\in U(1), g\in\mathcal{G}_E^1 ,
\]
where $c:\mathcal{G}_E^1\to\Z$ is the canonical continuous cocycle given by $c(x,k,y):=k$.
If $E$ is without sinks, then $\Psi_E$ restricts to an isomorphism
\[
C^*(E)^{U(1)}\cong
C^*(\mathcal{R}_E)
\]
where $\mathcal{R}_E$ is the topological subgroupoid of $\mathcal{G}_E$ with elements $(x,0,y)$, $x,y\in\partial E$, $x\sim_0 y$. Here we identify $C_c(\mathcal{R}_E^1)$ with functions in $C_c(\mathcal{G}_E^1)$ that are supported in $\mathcal{R}_E^1$. The above isomorphism is a special case of \cite[Lemma 3]{RRS17}, and the no sink condition is needed so that reduced and full groupoid C*-algebra are the same.

We shall write elements of $\mathcal{R}_E^1$ as pairs $(x,y)$, omitting the $0$ in the middle and think of this groupoid as a relation $\mathcal{R}_E^1\subseteq\partial E\times\partial E$.

\medskip

Let $E$ and $F$ be two directed graphs and assume that $\mathcal{G}_F$ is principal, i.e.~that the lag $k$ of a triple $(x,k,y)\in\mathcal{G}_F^1$ is uniquely determined by $x$ and $y$ (this happens for example when all vertices of $F$ are singular vertices). We can then omit the lag and think of 
$\mathcal{G}_F^1$ as a subset of $\partial F\times\partial F$.

\begin{rem}\label{rem:2.2}
Under the above assumptions, any isomorphism $\mathcal{G}_F\to\mathcal{R}_E$ has the form
\[
\mathcal{G}_F^1\ni (x,y)\mapsto\big(\omega(x),\omega(y)\big) \in\mathcal{R}_E^1 ,
\]
for some bijection $\omega:\partial F\to\partial E$. Conversely,
a bijection $\omega:\partial F\to\partial E$ induces:
\begin{enumerate}[label=(\arabic*)]
\item\label{rem:2.2en1} an isomorphism of algebraic groupoids if and only if $\omega\times\omega:\partial F\times \partial F\to \partial E\times\partial E$ maps the subset 
$\mathcal{G}_F^1$ to $\mathcal{R}_E^1$;

\item\label{rem:2.2en1} an isomorphism of topological groupoids if and only if, in addition, $\omega\times \omega$ (not just $\omega$!) is a homeomorphism.\vspace*{-16pt}
\end{enumerate}
\end{rem}

An isomorphism $\omega\times\omega:\mathcal{G}_F\to\mathcal{R}_E$ of topological groupoids induces an isomorphism $C^*(\mathcal{G}_F)\to C^*(\mathcal{R}_E)$ of full groupoid C*-algebras by defined on generators by
\begin{equation}\label{eq:iso}
1_X\longmapsto 1_{(\omega\times \omega)(X)} \quad\forall X\subseteq \mathcal{G}_F^1.
\end{equation}
In fact, the above isomorphism restricts to an isomorphism of Steinberg algebras \cite{R20}, which means that in the graph C*-algebra model this isomorphism is given by a polynomial formula in the Cuntz-Krieger generators.

\section{The groupoid of the quantum projective plane}\label{sec:5}

In this section, $F:=\widetilde{F}_2$ and $E:=\widetilde{L}_5$.
We label the vertices and edges of the graphs as in Section \ref{sec:3}, see \eqref{eq:shown}. Notice that $\partial F=F^\star$, since every vertex of $F$ is singular and there are no infinite paths, and in particular $\mathcal{G}_F$ is principal.
On the other hand, $\partial E=E^\infty$ since $E$ has no singular vertices, and we are in the situation described in Remark \ref{rem:2.2}.

We now exibith an isomorphism $\mathcal{G}_F\to\mathcal{R}_E$ of topological groupoids inducing the isomorphism of C*-algebras in Prop.~\ref{prop:isoB}.
As a byproduct, we get an independent proof of the results in Section~\ref{sec:3}.

\begin{thm}\label{thm:41}
Let $\omega:\partial F\to\partial E$ be the map
\begin{align*}
w_1 & \longmapsto e_1^\infty ,
&
w_2 & \longmapsto e_3^\infty ,
&
w_3 & \longmapsto e_5^\infty ,
\\
x_j & \longmapsto e_1^je_2e_3^\infty ,
&
y_k & \longmapsto e_3^ke_4e_5^\infty ,
&
x_jy_k & \longmapsto 
\begin{cases}
 e_1^{j+1}e_2e_3^ke_4e_5^\infty & \text{if } k<\binom{j+2}{2}, \\[2pt]
 e_1^je_2e_3^{k-j-1}e_4e_5^\infty & \text{if } k\geq \binom{j+2}{2}.
\end{cases}
\end{align*}
Then, $\omega\times\omega:\mathcal{G}^1_F\to\mathcal{R}^1_E$ is an isomorphism of topological groupoids.
\end{thm}

We split the proof in two steps. We first show that $\omega\times\omega:\mathcal{G}^1_F\to\mathcal{R}^1_E$ is an isomorphism of algebraic groupoids (Lemma \ref{lemma:bije}), and then that it is a homeomorphism (Lemma \ref{lemma:top}).

\begin{lemma}\label{lemma:bije}
\begin{enumerate}
\item\label{en:p1}
The map $\omega$ in Theorem \ref{thm:41} is bijective, with inverse given by the map
\begin{align*}
e_1^\infty & \longmapsto w_1 ,
&
e_3^\infty & \longmapsto w_2 ,
&
e_5^\infty & \longmapsto w_3 ,
\\
e_1^je_2e_3^\infty & \longmapsto x_j ,
&
e_3^ke_4e_5^\infty & \longmapsto y_k ,
&
e_1^je_2e_3^ke_4e_5^\infty & \longmapsto 
\begin{cases}
x_{j-1}y_k & \text{if } k<\binom{j+1}{2}, \\[2pt]
x_jy_{k+j+1} & \text{if } k\geq \binom{j+1}{2}.
\end{cases}
\end{align*}
\item\label{en:p2}
$\omega\times\omega:\mathcal{G}^1_F\to\mathcal{R}^1_E$
is an isomorphism of algebraic groupoids.
\end{enumerate}
\end{lemma}

\begin{proof}
Point \ref{en:p1} follows from Lemma \ref{lemma:sigma}.
Since $\omega$ maps paths in $F^\star$ with the same target into paths in $E^\infty$ with the same tail, from Remark \ref{rem:2.2}\ref{rem:2.2en1} it follows point \ref{en:p2}.
\end{proof}

\begin{lemma}\label{lemma:top}
$\omega\times\omega:\mathcal{G}^1_F\to\mathcal{R}^1_E$ is a homeomorphism.
\end{lemma}

\begin{proof}
If $\mathcal{G}$ and $\mathcal{R}$ are any two étale groupoids, and $\Omega:\mathcal{G}\to\mathcal{R}$ is an isomorphism of the underlying algebraic groupoids, then for all $U,V\subseteq\mathcal{G}$ one has
\[
\Omega(U\cdot V)=\Omega(U)\cdot\Omega(V) ,
\qquad
\Omega(U)^{-1}=\Omega(U)^{-1} .
\]
In particular, if $U,V$ are open in $\mathcal{G}$ and 
$\Omega(U),\Omega(V)$ are open in $\mathcal{R}$ (which is not automatic, since we are not assuming that $\Omega$ is an open map), since the multiplication in an étale groupoid is continuous,
the product $\Omega(U\cdot V)$ is also open in $\mathcal{R}$. To prove that $\Omega$ is an open map,
it is then enough to check that generators of the inverse semigroup of bisections have open image.
The same holds for $\Omega^{-1}$.

We now apply this argument to the map $\Omega:=\omega\times\omega$, with $\omega$ as in Lemma \ref{thm:41}.

Since $E$ has no singular vertices, a basis for the topology of $\mathcal{R}^1_E$ is given by the sets $Z_E(\mu,\nu)$ with $\mu,\nu\in E^\star$ any two finite paths with $t(\mu)=t(\nu)$ (by Remark \ref{rem:resim}\ref{rem:resimA}) and $|\mu|=|\nu|$ (by definition of the subgroupoid $\mathcal{R}_E$). Called $\ell$ the loop at $t(\mu)=t(\nu)$ and $k:=|\mu|=|\nu|$, then
\[
Z_E(\mu,\nu)=Z_E(\mu,\ell^k)\cdot Z_E(\nu,\ell^k)^{-1} .
\]
By the above argument, if we prove that $\Omega^{-1}(U)$ is a open for all $U$ of the form $Z_E(\mu,\ell^k)$, then it follows that $\Omega^{-1}(Z_E(\mu,\nu))$ is open for all $\mu,\nu$ and so $\Omega$ is continuous. We need to compute $\Omega^{-1}(Z_E(\mu,\ell^k))$
for all possible $\mu\in E^\star$. If the target of $\mu$ is $v_3$, then $Z_E(\mu,\ell^k)=\{(\mu e_5^\infty,e_5^\infty)\}$ is a singleton and its image under $\Omega^{-1}$ is a singleton $Z_F(\nu,w_3)=\{(\nu,w_3)\}$ for a suitable $\nu\in F^\star$ with target $w_3$ that can be computed using Lemma \ref{lemma:bije}. The remaining cases to be checked are ($p,q\geq 0$):
\begin{align*}
\mu &= e_1^{p} & \text{(source $v_1$ and target $v_1$)} , \\
\mu &= e_1^pe_2e_3^q & \text{(source $v_1$ and target $v_2$)} , \\
\mu &= e_3^{p} & \text{(source $v_2$ and target $v_2$)} .
\end{align*}
We have
\begin{multline*}
Z_E(e_1^p,e_1^p) =\{(e_1^\infty,e_1^\infty)\}\cup\{(e_1^je_2e_3^\infty,e_1^je_2e_3^\infty)\,|\,j\geq p\} \\
\cup
\big\{ (e_1^je_2e_3^ke_4e_5^\infty,e_1^je_2e_3^ke_4e_5^\infty)\,|\,j\geq p ,k\in\N \big\} ,
\end{multline*}
which gives
\begin{align*}
\Omega^{-1}(Z_E(e_1^p,e_1^p))
&=\{(w_1,w_1)\}\cup\{(x_j,x_j)\,|\,j\geq p\} \\
&\quad\cup
\big\{ (x_{p-1}y_k,x_{p-1}y_k)\,\big|\,j'\geq p-1 , k<\tbinom{j'+2}{2}\big\} \\
&\quad\cup
\big\{ (x_jy_k,x_jy_k)\,\big|\,j\geq p, k'\geq\tbinom{j+1}{2}+j+1=\tbinom{j+2}{2}\big\} ,
\end{align*}
where we renamed $j'=j-1$ and $k'=k+j+1$. 
We can regroup the elements and get
\begin{align*}
\Omega^{-1}(Z_E(e_1^p,e_1^p))
&=\{(w_1,w_1)\}\cup\{(x_j,x_j)\,|\,j\geq p\}\cup
\big\{ (x_jy_k,x_jy_k) \,\big|\,j\geq p, k\in\N\big\} \\
& \quad\cup\big\{ (x_{p-1}y_k,x_{p-1}y_k)\,\big|\, k<\tbinom{p+1}{2}\big\} \\
&=Z_F(w_1,w_1,S) \cup \bigcup_{k<\tbinom{p+1}{2}}Z_F(x_{p-1}y_k,x_{p-1}y_k) ,
\end{align*}
where $S=\{x_j\,|\,j<p\}$. Next, a similar computation gives:
\begin{align*}
\Omega^{-1}(Z_E(e_1^pe_2e_3^q,e_3^{p+q+1})) & = 
\{(x_p,w_2)\}\cup
\big\{ (x_{p-1}y_k,y_{p+k+1}) \,\big|\,  q\leq k<\tbinom{p+1}{2} \big\} \\
& \quad\cup
\big\{ (x_py_{p+k+1},y_{p+k+1}) \,\big|\,  k\geq q, k\geq\tbinom{p+1}{2} \big\} .
\end{align*}
Note that, if $q\geq\binom{p+1}{2}$ the second set in the above union is empty,
and in the third set the second side inequality is redundant. In this case, we get
\begin{equation}\label{eq:invertiloA}
\Omega^{-1}(Z_E(e_1^pe_2e_3^q,e_3^{p+q+1})) = Z_F(x_p,w_2,S) ,
\end{equation}
with $S=\{y_j\,|\,j<p+q+1\}$.
On the other hand, if $q<\binom{p+1}{2}$, then
\[
\Omega^{-1}(Z_E(e_1^pe_2e_3^q,e_3^{p+q+1}))  = 
Z_F(x_p,w_2,S)\cup\bigcup_{q\leq k<\tbinom{p+1}{2}}Z_F(x_{p-1}y_k,y_{p+k+1}) ,
\]
where $S=\{y_j\,|\,j<\binom{p+2}{2} \}$. Finally,
\begin{equation}\label{eq:invertiloB}
\Omega^{-1}(Z_E(e_3^{p},e_3^{p})) = Z_F(w_2,w_2,S) ,
\end{equation}
with $S=\{y_j\,|\,j<p\}$. This completes the proof that $\Omega$ is continuous.

We now prove that $\Omega$ is an open map. It follows from Remark \ref{rem:resim}\ref{rem:resimC} and from the above discussion that it is enough to check that
$\Omega(V)$ is a open for $V=Z_F(\mu,v,S)$, for all $\mu\in F^\star$ such that $v:=t(\mu)\in\{w_1,w_2\}$ is an infinite emitter, and for $V=Z_F(\mu,v)$, when $v=w_3$ is a sink. In the latter case,
$\Omega(Z_F(\mu,w_3))=\{(\nu,e_5^\infty)\}=Z_E(\nu e_5^\infty,e_5^\infty)$ for some $\nu\in E^\star$ with target $v_3$. It remains to compute the image of $Z_F(\mu,v,S)$ in the following cases
\begin{align*}
\mu &= w_1 & \text{(source $w_1$ and target $w_1$)} , \\
\mu &= x_j & \text{(source $w_1$ and target $w_2$)} , \\
\mu &= w_2 & \text{(source $w_2$ and target $w_2$)} .
\end{align*}
In the first case, it is enough to do the computation for $S=A_{k_0}:=\{x_k\,|\,k< k_0\}$ and for all $k_0\in\N$. In the other two cases, it is enough to do the computation for $S=B_{k_0}:=\{y_k\,|\,k< k_0\}$ and for all $k_0\in\N$. Indeed, every other open set $Z_F(\mu,v,S)$ with $S$ finite can be written as a union of a set of the above form and finitely many singletons of the form $Z_F(\nu,w_3)$. In fact, it is enough to do the check for all $k_0$ ``big enough''
(for all $k_0\geq N$ and for an arbitrary choice of $N$).

From \eqref{eq:invertiloA} and \eqref{eq:invertiloB} we get, for all $i\in\N$
and all $k_0\geq i$ (in this sense, ``big enough''), one has:
\begin{align*}
\Omega (Z_F(x_i,w_2,B_{k_0})) &=Z_E(e_1^ie_2^{k_0-i-1},e_3^{k_0}) ,
\\
\Omega (Z_F(w_2,w_2,B_{k_0})) &=Z_E(e^{k_0},e^{k_0}) .
\end{align*}
Finally, for all $p\in\N$,
\[
Z_F(w_1,w_1,A_p) =\{(w_1,w_1)\}\cup \{(x_j,x_j)\,|j\geq p\}\cup
\big\{ (x_jy_k,x_jy_k)\,\big|\,j\geq p,k\in\N \big\} .
\]
Then,
\begin{multline*}
\Omega(Z_F(w_1,w_1,A_p))=\{(e_1^\infty,e_1^\infty)\}\cup \{(e_1^je_2e_3^\infty,e_1^je_2e_3^\infty)\,|j\geq p\} \\
\cup
\big\{ (
e_1^{j'}e_2e_3^ke_4e_5^\infty,
e_1^{j'}e_2e_3^ke_4e_5^\infty
)\,\big|\,j'\geq p+1,k<\tbinom{j'+1}{2} \big\} \\
\cup
\big\{ (
e_1^je_2e_3^{k'}e_4e_5^\infty,
e_1^je_2e_3^{k'}e_4e_5^\infty
)\,\big|\,j\geq p,k'\geq\tbinom{j+1}{2} \big\} ,
\end{multline*}
where we renamed $j'=j+1$ and $k'=k-j-1$. We find:
\[
\Omega(Z_F(w_1,w_1,A_p))=Z_E(e_1^{p+1},e_1^{p+1})\cup Z_E\Big(e_1^pe_2e_3^{\binom{p+1}{2}},e_1^pe_2e_3^{\binom{p+1}{2}}\Big)
\]
if $p\geq 1$, and
\begin{equation}\label{eq:w1w1}
\Omega(Z_F(w_1,w_1))=Z_E(v_1,v_1)
\end{equation}
if $p=0$ (to prove \eqref{eq:w1w1} we observe that when $k<\binom{j'+1}{2}$, then $j'$ cannot be zero).
\end{proof}

\begin{thm}
The groupoid isomorphism $\Omega=\omega\times\omega$ in Theorem \ref{thm:41} induces the isomorphism of C*-algebras $\Phi$ of Section \ref{sec:3}.
\end{thm}

\begin{proof}
The isomorphism of groupoid C*-algebras is given by \eqref{eq:iso}. We now compute this map on generators. From \eqref{eq:w1w1} we see that
$1_{Z_F(w_1,w_1)}\longmapsto 1_{Z_E(v_1,v_1)}$. Next, $\Omega$ maps the singleton $Z_F(w_3,w_3)=\{(w_3,w_3)\}$ to 
$Z_E(v_3,v_3)=\{(e_5^\infty,e_5^\infty)\}$. Since the sum of vertex projection is $1$, we get
\[
1_{Z_F(w_i,w_i)}\longmapsto 1_{Z_E(v_i,v_i)} , \quad\forall\;i=1,2,3.
\]
Under the identification \eqref{eq:SZ}, this is exactly the first equality in \eqref{eq:elementsA}.
Next, the singleton $Z_F(y_j,w_3)=\{(y_j,w_3)\}$ is sent by $\Omega$ to the singleton
$Z_E(e_3^je_4,e_5^{j+1})=(e_3^je_4e_5^\infty,e_5^\infty)$. This gives the second equality in 
\eqref{eq:elementsA}. Finally,
\[
Z_F(x_j,w_2)=\{(x_j,w_2)\}\cup\{(x_jy_k,y_k)\,|\,k\in\N \} .
\]
Applying $\Omega$ we get
\begin{multline*}
\Omega(Z_F(x_j,w_2))=
\{(e^j_1e_2e_3^\infty,e_3^\infty)\}\cup
\big\{(e_1^je_2e_3^{k-j-1}e_4e_5^\infty,
e_3^ke_4e_5^\infty)\,|\,,k\geq\tbinom{j+2}{2} \big\} \\
\cup\big\{(
e_1^{j+1}e_2e_3^ke_4e_5^\infty
,e_3^ke_4e_5^\infty)\,|\,k<\tbinom{j+2}{2} \big\}  \\
=
Z_E(e_1^je_2e_3^{\binom{j+1}{2}},e_3^{\binom{j+2}{2}})
\cup \bigcup_{k<\binom{j+2}{2}}Z_E(e_1^{j+1}e_2e_3^ke_4
,e_3^ke_4e_5^{j+2}) .
\end{multline*}
Since the latter union is disjoint, its characteristic function is the sum of the characteristic functions of the individual sets, and we get
\begin{gather*}
  1_{Z_F(x_j,w_2)}\longmapsto 1_{Z_E(e_1^je_2e_3^{\binom{j+1}{2}},e_3^{\binom{j+2}{2}})}
    + \sum_{k=0}^{\binom{j+2}{2}-1} 1_{Z_E(e_1^{j+1}e_2\,e_3^ke_4,e_3^ke_4e_5^{j+2})} \;.
\end{gather*}
Under the identification \eqref{eq:SZ}, this exactly \eqref{eq:elementsB}.
\end{proof}

\section{An isomorphism of algebraic groupoids}\label{sec:6}

Let $R$ be a finite directed acyclic graph without multiple edges, $E$ the graph obtained from $R$ by adding a loop at every vertex, $F$ the graph obtained from $R$ by replacing every arrow with infinitely many,
$\mathcal{G}_F$ the boundary path groupoid of $F$ and $\mathcal{R}_E$ the subgroupoid of the boundary path groupoid of $E$ consisting of pairs of boundary paths that are tail equivalent with lag $0$. 
In formulas, one has
\[
E^0=F^0=R^0 , \qquad
F^1=R^1\times\N , \qquad
E^1=R^1\sqcup \{\ell_v\,|\,v\in R^0\} ,
\]
for all $e\in R^1,n\in \N$.
The edge $e\in R^1$ has the same source and target
in the graph $R$ and in the graph $E$,
the loop $\ell_v\in E^1$ has source and target $v\in E^0$, and, for all $n\in\N$ and $e\in E^1$, the edge $(e,n)\in F^1$ has the same source and target as $e$.
The aim of this section is to prove that
there is an isomorphism $\mathcal{G}_F\cong\mathcal{R}_E$ of algebraic groupoids.

\medskip

Note that $\partial E=E^\infty$ since $E$ has no singular vertices, and $\partial F=F^\star$ since there are no infinite paths in $F$ (it is finite and without directed cycles) and every vertex is singular. For every $(x,k,y)\in \mathcal{G}_F^1$, the lag $k=|x|-|y|$ is uniquely determined by $x$ and $y$. Thus $\mathcal{G}_F$ is principal, and we will write its elements as pairs $(x,y)$ omitting the lag.

We define a map $f:\partial F=F^\star\to\partial E=E^\infty$ as follows. If $\alpha=v$ is a vertex,
then $f(\alpha):=\ell_v^\infty$. If
\[
\alpha=(e_1,n_1)(e_2,n_2)\ldots (e_k,n_k)\in F^k , \quad k\geq 1,
\]
then
\begin{equation}
f(\alpha):=\ell_{s(e_1)}^{n_1}e_1\ell_{s(e_2)}^{n_2}e_2
\ldots\ell_{s(e_k)}^{n_k}e_k\ell_{t(e_k)}^\infty .
\label{eq:deff}
\end{equation}

\begin{prop}
The map $f\times f:\mathcal{G}^1_F\to\mathcal{R}^1_E$ is an isomorphism of algebraic groupoids.
\end{prop}

\begin{proof}
Two paths $\alpha,\beta\in F^\star$ are tail equivalent if and only if they have the same target. By construction,
if $t(\alpha)=t(\beta)$, then $f(\alpha)$ and $f(\beta)$ are tail equivalent with lag $0$ (they end with infinitely many loops at the same vertex). Thus, $f\times f$ maps $\mathcal{G}_F^1$ to $\mathcal{R}_E^1$.

The map $f:F^\star\to E^\infty$ is the composition of a graph homomorphism $F^\star\to E^\star$ and a suitable map $E^\star\to E^\infty$. The graph homomorphism is the identity on vertices, given edges by $(e,n)\mapsto \ell_{s(e)}^ne$, and extended to paths of positive length in the standard way. It is injective with image the set of all paths in $E^\star$ whose last edge is not a loop. The map $E^\star\to E^\infty$ adds a tail of the form $\ell_v^\infty$ to a path, and is a bijection between the set of paths whose last edge is not a loop, and the set $E^\infty$. The map $f$ is then bijective.

In view of Remark \ref{rem:2.2}, $f\times f$ is an isomorphism of algebraic groupoids.
\end{proof}

\begin{prop}\label{prop:homeo}
The map $f:\partial F\to\partial E$ is a homeomorphism.
\end{prop}

\begin{proof}
Since $f$ is a bijection between compact Hausdorff spaces, it is enough to prove that $f$ is continuous, i.e.\ that $f^{-1}$ is open. Since $E$ has no singular vertices, 
a basis of open sets for $\partial E$ consists solely of the cylinder sets 
\[
C_E(\alpha'):=\big\{\alpha' x\ell_{t(x)}^\infty\,\big|\,x\in E^\star,s(x)=t(\alpha')\big\} ,
\]
with $\alpha'\in E^\star$ (see Remark \ref{rem:similarly}).
Let
\[
\alpha'=\ell_{s(e_1)}^{n_1}e_1\ell_{s(e_2)}^{n_2}e_2
\ldots e_i\ell_{t(e_i)}^{n_{i+1}} ,
\]
where $e_1\ldots e_i$ is a path in $R$, $n_1,\ldots,n_{i+1}\in\N$ and $i\geq 1$.
Then,
the preimage under $f$ of $C_E(\alpha')$ is the set
$C_F(\alpha)\setminus\bigcup_{h\in S}C_F(\alpha h)$
where
\[
\alpha:=(e_1,n_1)\ldots (e_i,n_i)
\]
and
\[
S:=\big\{(x,k)\in R^1\times\N \,\big|\,s(x)=t(e_i),0\leq k< n_{i+1} \big\}.
\]
Finally, let $\alpha'=\ell_v^n$. Then,
the preimage under $f$ of $C_E(\alpha')$ is the set
$C_F(v)\setminus\bigcup_{h\in S}C_F(h)$
where now
\[
S:=\big\{(h,k)\in R^1\times\N \,\big|\,s(h)=v,0\leq k<n \big\}.
\]
Since the preimage of every basic open set is open, $f$ is continuous.
\end{proof}

Note that Prop.~\ref{prop:homeo} does not imply that $f\times f:\mathcal{G}^1_F\to\mathcal{R}^1_E$ is a homeomorphism (hence an isomorphism of topological groupoids) because the topology on $\mathcal{G}^1_F$ (resp.~$\mathcal{R}^1_E$) is not induced by the product topology on 
$\partial F\times\partial F$
(resp.~$\partial E\times\partial E$). In fact,

\begin{prop}
If the graph $R$ has a vertex that is neither a source nor a sink, then the map $f^{-1}\times f^{-1}:\mathcal{G}^1_E\to\mathcal{R}^1_F$ is not continuous.
\end{prop}

\begin{proof}
If $R^0$ has a vertex that is not a source nor a sink in $R$, this means that $E$ has a subgraph 
isomorphic to $\widetilde{L}_5$. We label the vertices and edges of this subgraph as in \eqref{eq:shown}.

Let $j\in\N$.
In $\mathcal{R}_E$, a basic open sets $Z_E(\mu,\nu)$ contains the point $(e_2e_3^\infty,e_3^\infty)$ if and only if 
$\mu$ is an initial subpath of $e_2e_3^\infty$ and $\nu$
is an initial subpath of $e_3^\infty$. Hence, $(\mu,\nu)=(e_2e_3^k,e_3^{k+1})$ and all basic open neighbourhoods of $(e_2e_3^\infty,e_3^\infty)$ have the form
\[
Z_E(e_2e_3^k,e_3^{k+1}) =\big\{(e_2e_3^kz,e_3^{k+1}z)\,\big|\,z\in E^\infty,s(z)=v_2 \big\}
\]
with $k\in\N$.
We deduce that
\begin{equation}\label{eq:seq}
(e_2e_3^ne_4e_5^\infty,e_3^{n+1}e_4e_5^\infty)
 \xrightarrow{\quad n\to\infty\quad} (e_2e_3^\infty,e_3^\infty) ,
\end{equation}
because the element on the left belongs to 
$Z_E(e_2e_3^k,e_3^{k+1})$ for all $n\geq k$, and then the sequence in \eqref{eq:seq} is eventually contained in any basic open neighbourhood of the point on the right hand side of \eqref{eq:seq}.

Let us call $x_n:=(e_2,n)$ and $y_n:=(e_4,n)$ the arrows in $F$.
Applying $f^{-1}\times f^{-1}$ to both sides of \eqref{eq:seq}, assuming that this function is continuous, we find that
the sequence $(x_0y_n,y_{n+1})_{n\in\N}$ converges to the point $(x_0,v_2)$.
Therefore, it must be eventually contained in the basic open neighbourhood
\[
Z_F(x_0,v_2)=\{ (x_0,v_2) \}\cup\big\{ (x_0y_k,y_k) \,\big|\,k\in \N\} .
\]
It follows that, for some $n$ big enough, there exists $k$ such that
\begin{equation}\label{eq:fix}
(x_0y_n,y_{n+1})=(x_0y_k,y_k) ,
\end{equation}
which is clearly impossible.
\end{proof}

It is not obvious how to modify the map $f$ in \eqref{eq:deff} in order to get an isomorphism of topological groupoids. If $\Omega=\omega\times\omega$ is such an isomorphism, then $f^{-1}\circ\omega=:\sigma:F^\star\to F^\star$ is a permutation of the (finite) paths in $F$ preserving the target. One should find a permutation that, composed with $f$, fixes the lack of continuity of $f^{-1}\times f^{-1}$. In the example $F:=\widetilde{F}_2$ and $E:=\widetilde{L}_5$ studied in Section \ref{sec:5}, if $\omega$ is the map in Theorem \ref{thm:41}, then $\sigma$ is the permutation that is the identity on paths of length $\leq 1$, and is given on paths of length $2$ by:
\[
x_jy_k \longmapsto 
\begin{cases}
 x_{j+1}y_k & \text{if } k<\binom{j+2}{2}, \\[2pt]
 x_jy_{k-j-1} & \text{if } k\geq \binom{j+2}{2}.
\end{cases}
\]
That is, it permutes the labels of length $2$ paths according to the  permutation of $\N^2$ in Lemma \ref{lemma:sigma}. Replacing $x_0y_n$ by $x_0y_{n-1}$ for large $n$ fixes the problem in \eqref{eq:fix}.

It is reasonable to believe that for $F:=\widetilde{F}_n$ and $E:=\widetilde{L}_{2n+1}$, $n\geq 2$, the map is that of the groupoids of $\widetilde{F}_{n-1}$
 and $\widetilde{L}_{2n-1}$ on all paths in $F$ of length $\leq n-1$, and something new happens only when considering paths in $F$ of length $n$. This problem will be investigated in future works.

\end{document}